\documentclass[11pt]{article}
\usepackage{amssymb, amstext, amscd, amsmath, color}

\usepackage{url}

\usepackage{amsfonts}
\usepackage{amsthm}
\usepackage{amsfonts}
\usepackage{array}

\usepackage{graphicx}
\usepackage{epstopdf}
\usepackage{epsfig}
\usepackage{eucal}
\usepackage{latexsym}
\usepackage{mathrsfs}
\usepackage{textcomp}
\usepackage{verbatim}
\usepackage{setspace}
\usepackage{subfigure}
\usepackage{tikz}
\usetikzlibrary{shapes.geometric}

\newtheorem{thm}{Theorem}[section]

\newtheorem{con}[thm]{Conjecture}

\newtheorem{defn}[thm]{Definition}

\numberwithin{equation}{section}

\newcommand{\pog}{\langle G, m \rangle} 
\newcommand{\pogp}{\langle G', \bm' \rangle}
\newcommand{\Tor}{\mathcal T_0} 
\newcommand{\pofw}{(\langle G, m \rangle,p)} 
\newcommand{\bm}{m} 
\newcommand{\gs}{\mathcal G_{\C}} 
\newcommand{\bbog}{\langle H, m \rangle} 

\newcommand{\bN}{{\mathbb{N}}}
\newcommand{\bQ}{{\mathbb{Q}}}
\newcommand{\bR}{{\mathbb{R}}}

  \newcommand{\C}{{\mathcal{C}}}

  \newcommand{\R}{{\mathcal{R}}}
\renewcommand{\S}{{\mathcal{S}}}
  \newcommand{\T}{{\mathcal{T}}}

\newcommand{\rank}{\operatorname{rank}}
\begin{document}

\title{One brick at a time: a survey of inductive constructions in rigidity theory}

\author{
{A. Nixon \thanks{tony.nixon@bristol.ac.uk, Heilbronn Institute for Mathematical Research, School of Mathematics, University of Bristol, U.K.} ~and E.  Ross \thanks{elissa@mathstat.yorku.ca, Department of Mathematics and Statistics, York University, Canada}
}
}

\date{}
\maketitle

\begin{abstract}
We present a survey of results concerning the use of inductive constructions to study the rigidity of frameworks. By inductive constructions we mean simple graph moves which can be shown to preserve the rigidity of the corresponding framework. We describe a number of cases in which characterisations of rigidity were proved by inductive constructions. That is, by identifying recursive operations that preserved rigidity and proving that these operations were sufficient to generate all such frameworks. We also outline the use of inductive constructions in some recent areas of particularly active interest, namely symmetric and periodic frameworks, frameworks on surfaces, and body-bar frameworks. As the survey progresses we describe the key open problems related to inductions. 
\end{abstract}

\section{Introduction}

Rigidity theory probes the question, given a geometric embedding of a graph, when is there a continuous motion or deformation of the vertices into a non-congruent embedding without breaking the connectivity of the graph or altering the edge lengths?
The geometric embeddings in question are typically bar-joint frameworks: collections of flexible joints and stiff bars that are permitted to pass through each other. 
The question of rigidity or flexibility is inherently dependent on the ambient space: in $1$ and $2$-dimensional Euclidean space there are complete combinatorial descriptions of the generic behaviour of a framework. In higher dimensions, however, there is no such characterisation; indeed there remains a number of challenging open problems.

In this survey we will concentrate on the most celebrated way of proving such a combinatorial description: an inductive construction. By an inductive construction we mean a constructive characterisation of a class of graphs or frameworks using simple operations. It is perhaps the simplicity of inductive constructions that make them so appealing, and helps to explain their widespread use. After all, the study of rigidity theory centres around highly intuitive concepts: building large rigid structures from smaller rigid components (e.g. building buildings from bricks). Inductive constructions provide an abstract analogue of this building-up process. 

There are two key ways in which inductive constructions have been used in rigidity theory. First, to show that a certain list of operations is sufficient to generate all graphs in a particular class (e.g. generic rigidity in the plane). Second, to show that certain inductive moves preserve rigidity (e.g. vertex splitting). As a result, inductive constructions have been used as proof techniques without necessarily hoping for complete combinatorial characterisations (e.g. the proof of the Molecular Conjecture). When a complete combinatorial description is obtained, inductive characterisations typically do not make for fast algorithms. On the other hand, once we have an inductive sequence for a rigid framework, we have an instant certificate of its rigidity.  

We begin the survey with a gentle introduction into rigidity and global rigidity theory in $2$-dimensions from an inductive perspective. From there we outline the key open problems in extending inductive constructions to $3$-dimensional frameworks before describing some purely graph theoretical inductive constructions in Section \ref{sec:klTight}.
The central topic of discussion in Sections \ref{sec:periodic} and \ref{sec:symmetric} is the rigidity of periodic and symmetric frameworks, two types of frameworks with special geometric features. Following that we discuss frameworks on surfaces and body-bar frameworks (Sections \ref{sec:surfaces} and \ref{sec:bodyBar}) before finishing the survey by briefly outlining, in Section \ref{sec:further}, a number of other avenues of rigidity theory which have benefitted from inductive techniques.

\section{Basics of Rigidity}

A {\it (bar-joint) framework} is an ordered pair $(G,p)$ where $G$ is a graph and $p:V\rightarrow \bR^d$ is an embedding of the vertices into $\bR^d$.
We are interested in the typical behaviour of frameworks. Thus we say that a framework is {\it generic} if the coordinates of the framework points form an algebraically independent set (over $\bQ$).
Two frameworks on the same graph $(G,p)$ and $(G,q)$ are {\it equivalent} if the (Euclidean) edge lengths in $(G,p)$ are the same as those in $(G,q)$ and are \emph{congruent} if the distance between pairs of points in $(G,p)$ are the same as those in $(G,q)$. 

\begin{defn}
A framework $(G,p)$ is {\rm flexible}  in $\bR^d$ if there is a continuous motion $x(t)$ of the framework points such that $(G,x(t))$ is equivalent to $(G,p)$ for all $t$ but is not congruent to $(G,p)$ for some $t$ (where $x(t)\neq p$). $(G,p)$ is {\rm (continuously) rigid} if it is not flexible.
\end{defn}

Understanding rigidity becomes more tractable after linearising the problem.
The \emph{rigidity matrix} $R_d(G,p)$ is a sparse matrix where each row corresponds to an edge, and with the appropriate ordering each $d$-tuple of columns corresponds to the coordinates of a framework vertex. The entries in row $ij$ are zero except in the columns corresponding to $i$ and $j$ where the entries are $p_i-p_j$ and $p_j-p_i$ respectively. This matrix is (up to scaling) the Jacobean derivative matrix of the system of quadratic edge length equations.
The (infinitesimal) rigidity matroid $\R_d$ (for a generic framework $(G,p)$) is the linear matroid induced by linear independence in the rows of the rigidity matrix $R_d(G,p)$.

\begin{defn}
Let $p=(p_1,\dots,p_{|V|})$. An infinitesimal flex $u=(u_1,\dots,u_{|V|}) \in \bR^{d|V|}$ is a vector satisfying $(p_i-p_j).(u_i-u_j)=0$ for all edges $ij$.
A framework is \emph{infinitesimally rigid} if there are no non-trivial infinitesimal flexes.
\end{defn}

There are examples of frameworks that are infinitesimally flexible but continuously rigid, however all such examples occur for geometric reasons.

\begin{thm}[Asimow and Roth \cite{AsimowRoth}]
Let $(G,p)$ be a generic framework. Then $(G,p)$ is (continuously) rigid if and only if it is infinitesimally rigid.
\end{thm}

\subsection{Constructing Frameworks in the Plane}

Let us, for now, restrict attention to frameworks in the plane and consider the following construction moves \cite{Henneberg}:

\begin{enumerate}
\item[1] add a vertex $v$ with $d(v)=2$ and $N(v)=\{a,b\}$, $a \neq b$, 
\item[2] remove an edge $xy$, $x\neq y$, and add a vertex $v$ with $d(v)=3$ and $N(v)=\{x,y,z\}$ for some $z \in V$,
\end{enumerate}
In the literature, operation 1 may be referred to as a {\it Henneberg 1} move, a {\it $0$-extension} or a {\it vertex addition} and operation 2 may be referred to as a {\it Henneberg 2} move \cite{PseudoTriangles}, \cite{NixonOwenPower}, \cite{Tayhenneberg}, a {\it $1$-extension} \cite{GraverServatius2}, \cite{JacksonJordan}, \cite{JordanKaszanitskyTanigawa} or {\it edge splitting} \cite{RossThesis}, \cite{Schulze}, \cite{Whiteleyscene}. All have their advantages: Henneberg move gives credit to the original source; $i$-extension indicates the number of edges removed in the operation; vertex addition and edge splitting are the most accessible to newcomers to the subject. Since this is a survey we choose to use the last option from here on.

\begin{figure}[ht]
\centering
\subfigure{
\includegraphics[width=2in]{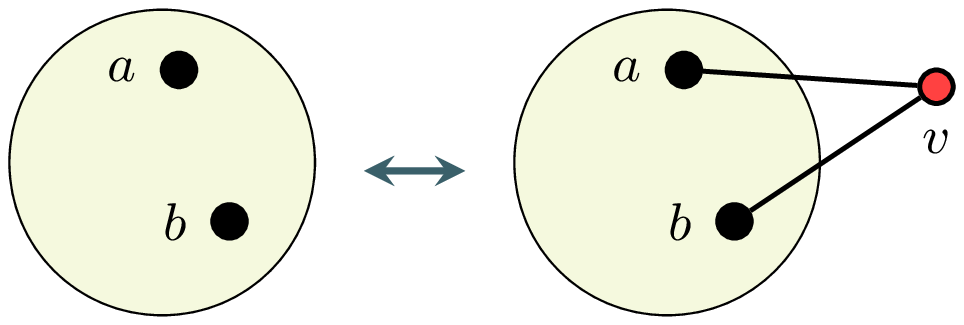}}
\hspace{.5in}
\subfigure{
\includegraphics[width=2in]{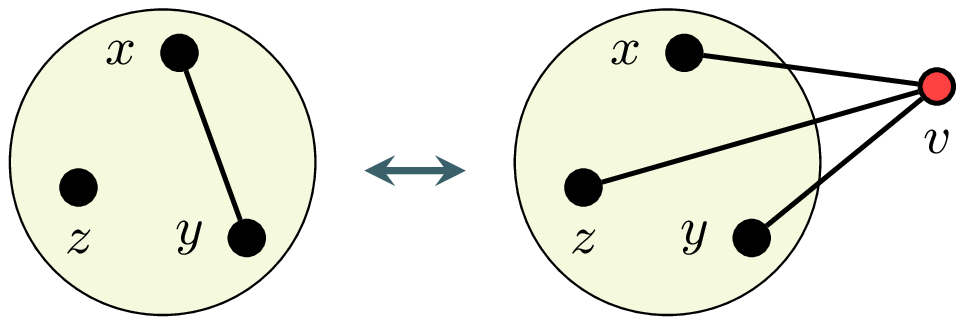}}
\caption{Vertex addition and Edge splitting.}

\end{figure}

\begin{defn}\label{def:23tight}
A graph $G=(V,E)$ is {\rm $(2,3)$-sparse} if for every subgraph $G'=(V',E')$ with at least one edge, $|E'|\leq 2|V'|-3$. $G$ is {\rm $(2,3)$-tight} if $G$ is $(2,3)$-sparse and $|E|=2|V|-3$.
\end{defn}

\begin{thm}[Henneberg \cite{Henneberg}, Laman \cite{Laman}]\label{thm:Hen}
A graph $G$ is $(2,3)$-tight if and only if it can be derived recursively from $K_2$ (the single edge) by vertex additions and edge splitting.
\end{thm}

\begin{figure}[ht]
\centering
\includegraphics[width=4in]{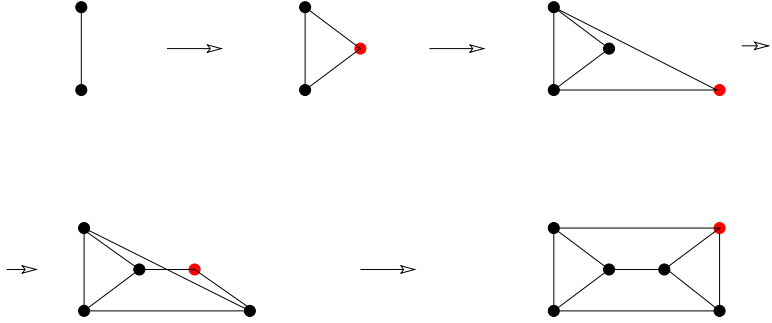}
\caption{A Henneberg-Laman sequence for the triangular prism.}
\end{figure}

An infinitesimally rigid graph $G$ is called {\it isostatic} or {\it minimally rigid} if deleting any edge will destroy its rigidity. 

\begin{thm}[Laman's theorem \cite{Laman}]
A graph $G$ is generically minimally rigid in the plane if and only if $G$ is $(2,3)$-tight.
\end{thm}

Maxwell \cite{Maxwell} proved that any generically minimally rigid graph must be $(2,3)$-tight. The harder sufficiency direction relies on Theorem \ref{thm:Hen}. Given the inductive construction, and since $K_2$ clearly has a generically rigid realisation, it remains only to show that the result of applying vertex addition and edge splitting to a generically minimally rigid graph is a generically minimally rigid graph. Vertex addition is trivial; we have a rigidity matrix with rank equal to the number of rows and genericness ensures the two rows and two columns, that we add, increases the rank by two. Edge splitting is slightly more involved. Let $G'$ be formed from $G$ by edge splitting. Then a typical proof uses the fact that for a graph $H=(V,E)$ and two maps $q_1, q_2:V\rightarrow \bR^2$ with $(H,q_1)$ generic, $\rank R_2(H,q_1) \geq \rank R_2(H,q_2)$ see, for example, \cite{Whiteleymatroid}. 
Using this, choose the new vertex in $G'$ to be on the line through the (not yet) removed edge in $(G',p')$. The collinear triangle created corresponds to a minimal set of linearly dependent rows in the rigidity matrix (i.e. a circuit in the rigidity matroid). We can remove any edge from this triangle without reducing the rank of $R_2(G+v,p')$. Since it is clear that $\rank R_2(G,p)+2=\rank R_2(G+v, p')$, it follows that $\rank R_2(G',p')=\rank R_2(G,p)+2$. 

\section{Global Rigidity}

There are a number of applications in which rigidity is not strong enough due to the possibility of multiple distinct realisations with the same edge lengths. Global rigidity corresponds exactly to there being a unique realisation, up to congruence, with the given edge lengths.
For a full survey on global rigidity see \cite{JacksonJordansurvey}, we give only a brief description of the use of inductive constructions for global rigidity. 

\begin{defn}
A framework $(G,p)$ is {\rm generically globally rigid} if for all equivalent choices of $q$ the frameworks $(G,p)$ and $(G,q)$ are congruent.
\end{defn}

In characterising global rigidity we will also use the following strong form of rigidity.

\begin{defn}
Let $G=(V,E)$. A framework $(G,p)$ is {\rm redundantly rigid} if $(G,p)$ is rigid and for all $e\in E$ the framework $(G-e,p)$ is rigid.
\end{defn}

\begin{thm}[Hendrickson \cite{Hendrickson}]\label{thm:hennecessary}
Let $(G,p)$ be a generic globally rigid framework in $\bR^d$. Then $G$ is a complete graph on at most $d+1$ vertices or $G$ is $(d+1)$-connected and $(G,p)$ is redundantly rigid in $\bR^d$.
\end{thm}

\subsection{Circuits}

By Laman's theorem the minimal number of edges needed for a graph to be generically globally rigid in the plane is $2|V|-2$. By Theorem \ref{thm:hennecessary} the graph must also be redundantly rigid. This implies that if $G$ is generically globally rigid with $2|V|-2$ edges then $G$ is a $(2,3)$-circuit; that is a graph with $2|V|-2$ edges in which every proper subgraph (with at least one edge) is $(2,3)$-sparse. 

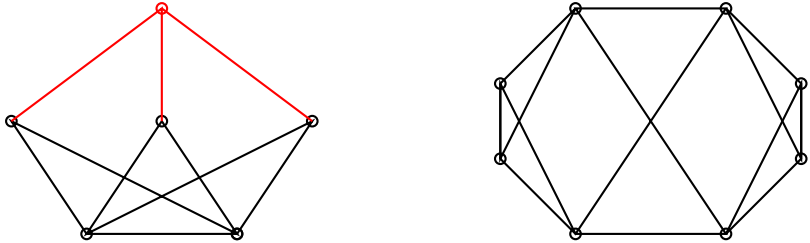
\begin{figure}[ht]
\begin{center}
\begin{tikzpicture}
\draw[black,thick] (0,0) circle (2pt);
\draw[black,thick] (-2,0) circle (2pt);
\draw[black,thick] (-3,1.5) circle (2pt);
\draw[black,thick] (-1,1.5) circle (2pt);
\draw[black,thick] (1,1.5) circle (2pt);
\draw[red,thick] (-1,3) circle (2pt);

 \draw[black,thick]
  (0,0) -- (-2,0) -- (-3,1.5) -- (0,0) -- (-1,1.5) -- (-2,0) -- (1,1.5) -- (0,0);
 \draw[red,thick]
 (-3,1.5) -- (-1,3) -- (-1,1.5);
 \draw[red,thick]
 (1,1.5) -- (-1,3);

\draw[black,thick] (4.5,0) circle (2pt);
\draw[black,thick] (6.5,0) circle (2pt);
\draw[black,thick] (4.5,3) circle (2pt);
\draw[black,thick] (6.5,3) circle (2pt);
\draw[black,thick] (0,0) circle (2pt);
\draw[black,thick] (3.5,2) circle (2pt);
\draw[black,thick] (3.5,1) circle (2pt);
\draw[black,thick] (7.5,2) circle (2pt);
\draw[black,thick] (7.5,1) circle (2pt);

\draw[black,thick]
  (4.5,0) -- (6.5,0) -- (4.5,3) -- (6.5,3) -- (4.5,0) -- (3.5,1) -- (3.5,2) -- (4.5,0);
\draw[black,thick]
 (4.5,3)-- (3.5,1) -- (3.5,2) -- (4.5,3);
\draw[black,thick]
 (6.5,3)-- (7.5,1) -- (7.5,2) -- (6.5,3);
\draw[black,thick]
 (6.5,0)-- (7.5,1) -- (7.5,2) -- (6.5,0);
\end{tikzpicture}
\end{center}
\caption{Two examples of $(2,3)$-circuits. On the left the red vertex cannot be reduced as the result will be a copy of $K_4$ with a degree 2 vertex adjoined. On the right there is no inverse edge splitting move that results in a $(2,3)$-circuit.}
\label{fig:circuit}
\end{figure}

The beauty of Theorem \ref{thm:Hen} is that for every $(2,3)$-tight graph containing a vertex of degree $3$, there is always an inverse edge splitting operation resulting in a smaller $(2,3)$-tight graph. This is not the case for $(2,3)$-circuits; hence it is significantly more challenging to prove an inductive construction.
For example it is possible for a degree $3$ vertex $v$ in a $(2,3)$-circuit to have all neighbours $x,y,z$ of degree $3$. Here any inverse edge splitting operation results in a graph with $2|V|-2$ edges which is not a circuit since at least one of $x,y,z$ has degree $2$, see Figure \ref{fig:circuit}.

\begin{thm}[Berg and Jord{\'a}n \cite{BergJordan}]
Let $G$ be a $3$-connected $(2,3)$-circuit. Then there is an inverse edge splitting move on some vertex of $G$ that results in a smaller $(2,3)$-circuit.
\end{thm}

Combining this with the well known 2-sum operation from matroid theory allowed them to inductively characterise $(2,3)$-circuits. The 2-sum operation glues two $(2,3)$-circuits together along an edge and deletes the common edge. The inverse operation separates along a $2$-vertex cutset. This operation has further been examined from the rigidity perspective in \cite{Servatius2sum}.

\begin{thm}[Berg and Jord{\'a}n \cite{BergJordan}]\label{thm:BergJordancircuit}
A graph $G$ is a $(2,3)$-circuit if and only if $G$ can be generated from copies of $K_4$ by applying edge splitting moves within connected components and taking 2-sums of connected components.
\end{thm}

While it was easy to see that the edge splitting operation preserves rigidity, showing that it preserves global rigidity is more intricate. This was originally proved by Connelly \cite{Connelly} as a corollary to his sufficient condition for global rigidity in terms of the rank of the stress matrix. An alternative proof was later given by Jackson, Jord{\'a}n and Szabadka \cite{JacksonJordanSzabadka} during their analysis of globally linked vertices.

\subsection{Characterising Global Rigidity}

The characterisation of $(2,3)$-circuits was extended by Jackson and Jord{\'a}n to $M$-connected graphs; these are graphs in which there is a $(2,3)$-circuit containing any pair of edges, i.e. the rigidity matroid is connected.
They showed using ear decompositions that all $3$-connected, $M$-connected graphs could be generated from $K_4$ by edge splitting operations and edge additions. Part of the subtlety here is that they had to be able to alternate between the operations, see \cite[Figure $6$]{JacksonJordan}.

\begin{thm}[Hendrickson \cite{Hendrickson}, Connelly \cite{Connelly}, Jackson and Jord{\'a}n \cite{JacksonJordan}]\label{thm:2dglobal}
A framework $(G,p)$ is generically globally rigid in the plane if and only if $G$ is a complete graph on at most 3 vertices or $G$ is $3$-connected and $(G,p)$ is redundantly rigid.
\end{thm}

\section{Rigidity in $3$-space}
\label{sec:3d}

\begin{center}
\begin{figure}[ht]
\centering
\includegraphics[height=3.5cm]{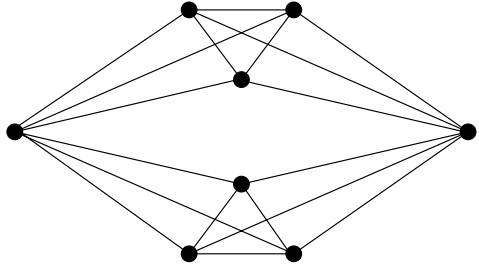}
\caption{The double banana; a flexible circuit in the $3$-dimensional rigidity matroid.}
\label{fig:double banana}
\end{figure}
\end{center}

As in the plane the necessity of combinatorial counts for minimal rigidity was shown by Maxwell \cite{Maxwell}. The appropriate graphs are the $(3,6)$-tight graphs. However it is no longer true that these graphs are sufficient for minimal rigidity; there exist $(3,6)$-tight graphs which are generically flexible in $3$-dimensions, see Figure \ref{fig:double banana} for an example. Thus the outstanding open problem in rigidity theory is to find a good combinatorial description of generic minimal rigidity in $3$-dimensions.

\vspace{-1cm}

\begin{center}
\begin{figure}[ht]
\centering
\includegraphics[width=2in]{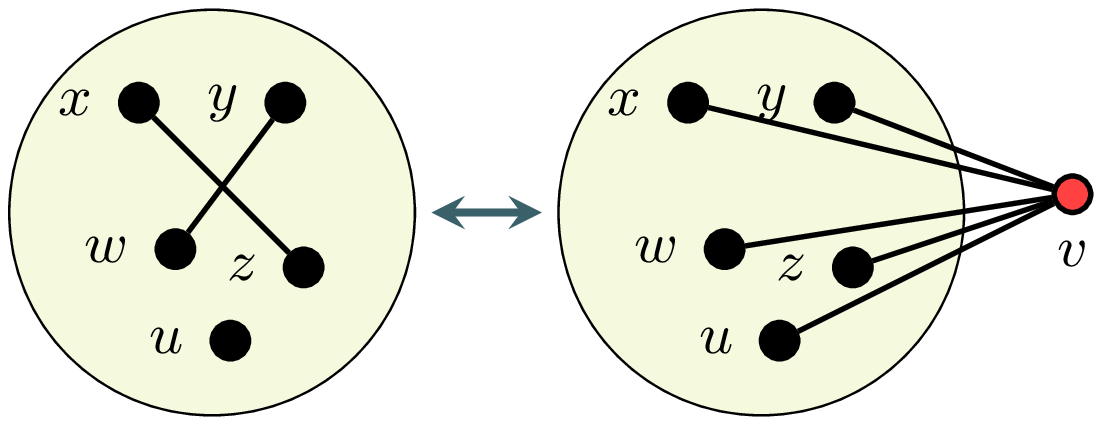}\hspace{.5in}
\includegraphics[width=2in]{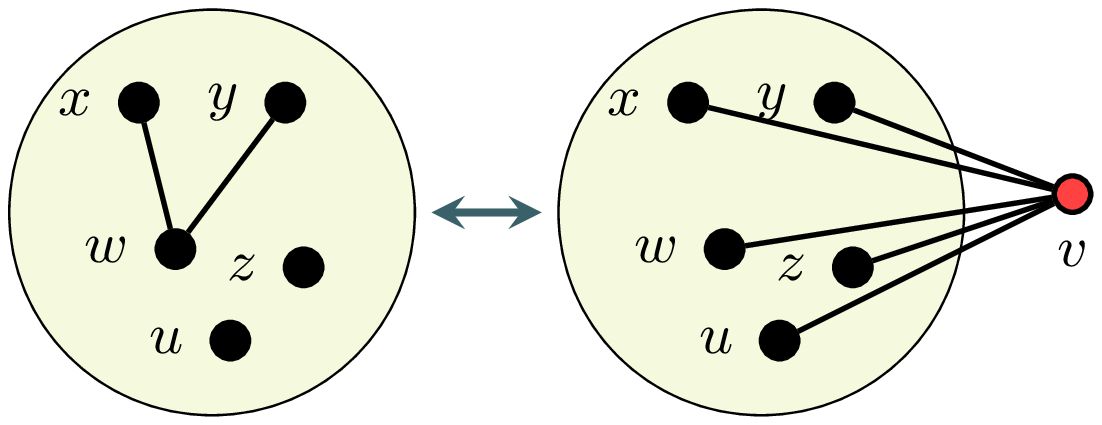}
\caption{The $X$- and $V$-replacement operations in three dimensions.}
\label{fig:xreplace}
\end{figure}
\end{center}

From an inductive construction perspective it is known that the analogues of vertex addition and edge splitting preserve rigidity. In fact Tay and Whiteley \cite{TayWhiteley} proved that, in dimension $d$, the addition of a vertex of degree $d$ (vertex addition) or the subdivision of an edge combined with adding $d-1$ additional edges incident to the new vertex (edge splitting) preserves rigidity.
However the average degree in a $(3,6)$-tight graph approaches $6$. Thus we require new operations to deal with degree $5$ vertices, see Figure \ref{fig:xreplace}.

\subsection{Degree $5$ Operations}

\vspace{-1cm}

\begin{con}[Whiteley \cite{Whiteleyscene}]
Let $G$ be generically rigid in $\bR^3$ and let $G'$ be the result of an $X$-replacement applied to $G$. Then $G'$ is generically rigid in $\bR^3$.
\end{con}

The conjecture is intuitively appealing since for the variant in the plane, similarly to the edge splitting argument, it is easy to establish the preservation of rigidity. 
Let $G'$ be formed from $G$ by an $X$-replacement. Then for some pair of edges in $G$, say $uv$ and $xy$, $G'$ is formed from $G-\{uv,xy\}$ by adding a single vertex $z$ and edges $uz,vz,xz,yz$. We choose a realisation $p^*$ of $G+z$ such that $z$ lies on the unique point defining the intersection of the lines through $uv$ and $xy$ (since we may assume $G$ was generic these lines are not parallel). Now in the rigidity matroid for $(G+z,p^*)$ there are two collinear circuits, defined by the edge sets $\{uv,uz,vz\}$ and $\{xy,xz,yz\}$ respectively. Thus the deletion of $uv$ and $xy$ does not reduce the rank of the rigidity matrix. The statement follows since it is not hard to argue that $\rank R_2(G+z,p^*)$ is 2 more than $\rank R_2(G,p)$.

However this argument easily breaks down in higher dimensions; generically two lines do not intersect.
Going against the conjecture are the following two facts; the analogue of $X$-replacement in $4$-dimensions fails and $X$-replacement in $3$-dimensions does not preserve global rigidity. 

The first fact is based on a general argument \cite{GraverServatius2}, which in particular shows that $K_{6,6}$ is dependent in the $4$-dimensional rigidity matroid.

\begin{center}
\begin{figure}[ht]
\centering
\includegraphics[width=6.5cm]{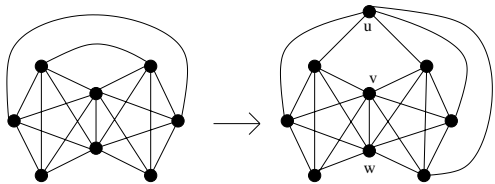}
\caption{An example due to Tibor Jord{\'a}n showing $X$-replacement does not necessarily preserve global rigidity in $3$-dimensions, \cite{Banff08}.}
\label{globalxreplace}
\end{figure}
\end{center}

The second fact is illustrated in Figure \ref{globalxreplace}. The first graph is generically globally rigid in $3$-dimensions. This is easily seen since it can be formed from $K_5$ by a sequence of ($3$-dimensional) edge splitting moves and edge additions, both of which preserve global rigidity. The second graph, obtained by an $X$-replacement on the first graph, contains a $3$-vertex-cutset $\{u,v,w\}$. Thus Theorem \ref{thm:hennecessary} implies it is not globally rigid.

\vspace{-1cm}

We also mention that \cite{GraverServatius2} gives a nice discussion of the problem including several special cases where $X$-replacement is known to preserve rigidity.

It is quickly apparent that the $V$-replacement operation presents a new difficulty; the earlier inductive operations were easily seen to preserve the relevant vertex/edge counts on the graph and all subgraphs. It is not true, however, that $V$-replacement always preserves the subgraph counts, we may make a bad choice of vertex $w$. Tay and Whiteley \cite{TayWhiteley} have made a double-$V$ conjecture but this has an immediate problem from an algorithmic perspective. The conjecture implies an inductive construction of minimally rigid graphs in $3$-space. Using this inductive construction to check the rigidity of a given framework would require recording both graphs each time the $V$-replacement is applied. Thus for worst case graphs the generating sequence of inductive operations requires remembering exponentially many different graphs.

\subsection{Vertex Splitting}

\begin{center}
\begin{figure}[ht]
\centering
\includegraphics[width=3.5in]{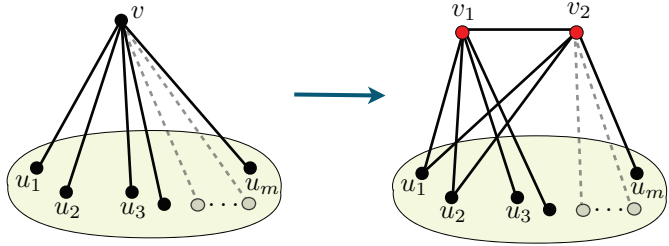}
\caption{The $3$-dimensional vertex splitting operation.}
\label{Vertex splitting}
\end{figure}
\end{center}

Let $v \in V$ have $N(v)=\{u_1, \dots, u_{m}\}$. A \emph{vertex splitting} operation (in $3$-dimensions) on $v$ removes $v$ and its incident edges, adds vertices $v_0,v_1$ and edges $u_1v_0,u_2v_0, u_1v_1,u_2v_1, v_0v_1$ and re-arranges the edges $u_3v,\dots,u_{m}v$ in some way into edges $u_iv_j$ for $i\in \{3,\dots,m\}$ and $j \in \{0,1\}$. See Figure \ref{Vertex splitting} and also \cite{Whiteleyvertexsplitting} where the operation was introduced for $d$-dimensional frameworks.

\begin{thm}[Whiteley \cite{Whiteleyvertexsplitting}]\label{thm:vertexsplittingd}
Let $G$ have a generically minimally rigid realisation in $\bR^d$ and let $G'$ be formed from $G$ by a vertex splitting operation. Then $G'$ has a generically minimally rigid realisation in $\bR^d$.
\end{thm}

For a globally rigid graph in the plane it can be derived from Theorem \ref{thm:2dglobal}, see \cite{JordanSzabadka}, that applying a vertex splitting operation, in which each new vertex is at least $3$-valent, results in a globally rigid graph.

\begin{con}[Cheung and Whiteley \cite{CheungWhiteley}]
Let $G$ be globally rigid in $\bR^d$ and let $G'$ be formed from $G$ by a vertex splitting operation such that each new vertex is at least $d+1$-valent. Then $G'$ is globally rigid in $\bR^d$.
\end{con}

Vertex splitting has also been used to prove a variety of results for restricted classes of three-dimensional frameworks. In particular, Finbow and Whiteley recently used vertex splitting to prove that {\it block and hole frameworks} are isostatic \cite{FinbowWhiteley}. A block and hole framework is a triangulated sphere (known to be isostatic by early results of Cauchy and Dehn) where some edges have been removed to create {\it holes}, while others added to create isostatic subframeworks called {\it blocks}, all the while maintaining the general $|E| = 3|V| - 6$ count. An example of such a framework is a geodesic dome. The base of the dome can be considered as a block. It becomes possible to remove some edges from the rest of the dome, perhaps to create windows and doors. The result of Finbow and Whiteley will identify which edges may be removed. The proof of this result relies on vertex splitting in a central way.

\section{Inductive constructions for $(k,l)$-tight graphs}
\label{sec:klTight}

Up until now it has been obvious that we concentrated on simple graphs i.e. graphs with no loops or multiple edges. From here on graphs will allow loops or multiple edges and we will specify that graphs which do not are simple.

\begin{defn}
Let $k,l \in \bN$ and $l < 2k$. A graph $G=(V,E)$ is {\rm $(k,l)$-sparse} if for every subgraph $G'=(V',E')$, with $|V'|\geq k$, $|E'|\leq k|V'|-l$. $G$ is {\rm $(k,l)$-tight} if $G$ is $(k,l)$-sparse and $|E|=k|V|-l$.
\end{defn}

We choose to restrict to the range $l<2k$ since in this range $(k,l)$-tight graphs are the bases of matroids \cite{Whiteleymatroid} and \cite{LeeStreinu}. Observe that $(3,6)$-tight graphs are outside this range; indeed they do not form the bases of a matroid.
Since we now allow multiple edges there are more possibilities for vertex additions and edge splitting operations. Throughout the rest of the paper, when we consider graphs these operations will be understood to allow the graph variants, see Figures \ref{fig:vertexAddition} and \ref{fig:edgeSplit}.

In \cite{FrankSzego}, Frank and Szeg\H{o} prove inductive characterisations of graphs which are {\it nearly $k$-tree connected}, which naturally extend the combinatorial elements of Henneberg's original result. 

\begin{defn}
A graph $G$ is called {\rm $k$-tree connected} if it contains $k$ edge-disjoint spanning trees. A graph is {\rm nearly $k$-tree connected} if it is not $k$-tree connected, but the addition of any edge to $G$ results in a $k$-tree connected graph. 
\end{defn}

Note that Henneberg's result (Theorem \ref{thm:Hen}) can be rephrased as follows: A graph $G$ is nearly $2$-tree connected if and only if it can be constructed from a single edge by a sequence of vertex additions and edge splitting operations.
\begin{thm}[Frank and Szeg\H{o} \cite{FrankSzego}]\label{thm:frankszego}
A graph $G$ is nearly $k$-tree-connected if and only if $G$ can be constructed from the graph consisting of two vertices and $k-1$ parallel edges by applying the following operations:
\begin{enumerate}\item add a new vertex $z$ and $k$ new edges ending at $z$ so that there are no $k$ parallel edges,\item choose a subset $F$ of $i$ existing edges ($1 \leq i \leq k - 1$), pinch the elements of $F$ with a new vertex $z$, and add $k - i$ new edges connecting $z$ with other vertices so that there are no $k$ parallel edges in the resulting graph.
\end{enumerate}
\end{thm}

We recall a result of Nash-Williams \cite{NashWilliams}, which states that a graph $G=(V,E)$ is the union of $k$ edge-disjoint forests if and only if $|E'| \leq k|V'| - k$ for all nonempty subgraphs $G' = (V', E') \subseteq G$. Continuing the theme of extending Henneberg's theorem, by this result of Nash-Williams, Frank and Szeg\H{o} show that a graph $G$ is $(k, k+1)$-tight if and only if it is nearly $k$-tree connected. We point the interested reader to a book of Recski \cite{Recskibook} where a number of connections between minimally rigid graphs and tree decompositions are proved.

Fekete and Szeg\H{o} have established a Henneberg-type characterisation theorem of $(k,l)$-sparse graphs for the range $0 \leq l \leq k$.  The following definition extends vertex addition and edge splitting to arbitrary dimension. 
\begin{defn}
Let $G$ be a graph, and let $0 \leq j \leq m \leq k$. Choose $j$ edges of $G$ and pinch into a new vertex $z$. Put $m-j$ loops on $z$, and link it with other existing vertices of $G$ by $k-m$ new edges. This move is called an {\rm edge pinch}, and will be denoted $K(k,m,j)$. 
\end{defn}
The graph on a single vertex with $l$ loops will be denoted $P_{l}$. The main result of \cite{FeketeSzego} is the following.
\begin{thm}[Fekete and Szeg\H{o} \cite{FeketeSzego}]
Let $G=(V,E)$ be a graph and let $1\leq l \leq k$. Then $G$ is a $(k,l)$-tight graph if and only if $G$ can be constructed from $P_{k-l}$ with operations $K(k,m,j)$ where $j \leq m \leq k-1$ and $m-j \leq k-l$. 

$G$ is a $(k,0)$-tight graph if and only if $G$ can be constructed from $P_k$ with operations $K(k,m,j)$, where $j \leq m \leq k$ and $m-j \leq k$. 
\end{thm}

This result has subsequently been applied to periodic body-bar frameworks \cite{RossBodyBar}, see Section \ref{sec:rossBodyBar}. Inductive moves for $(k,l)$-tight graphs have also been considered using an algorithmic perspective in \cite{LeeStreinu}.

\section{Periodic Frameworks}
\label{sec:periodic}
Over the past decade, the topic of periodic frameworks has witnessed a surge of interest in the rigidity theory community \cite{RoyalSociety,BorceaStreinu,BorceaStreinuII,MalesteinTheran,RossThesis}, in part due to questions raised about the structural properties of zeolites, a type of crystalline material with numerous practical applications. Inductive constructions have been used to provide combinatorial characterisations of certain restricted classes of periodic frameworks, which we describe below. 

\begin{figure}[h!]
\begin{center}
\subfigure{\includegraphics[width=1.5in]{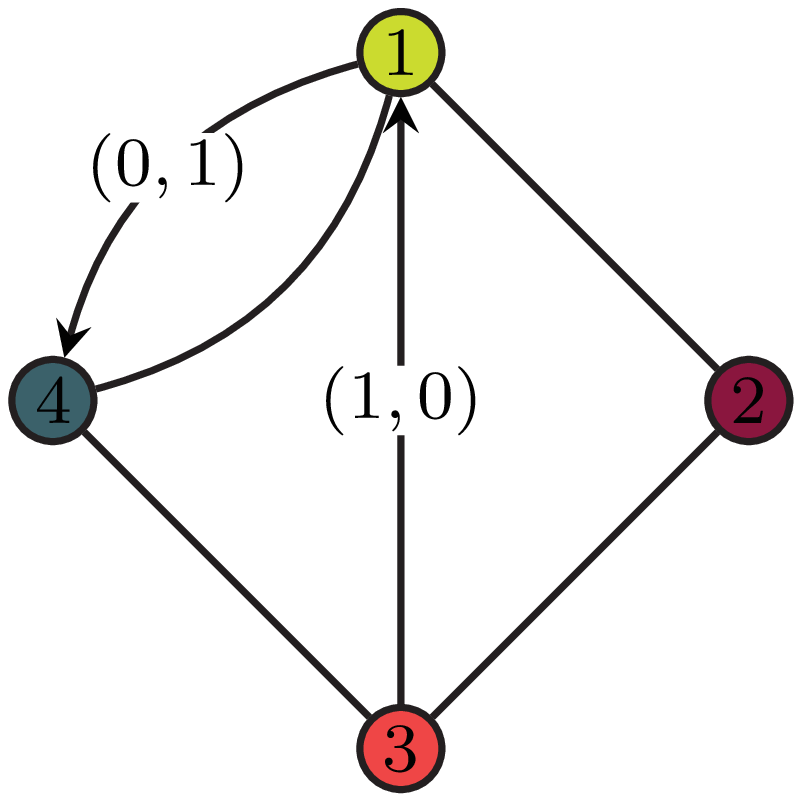}}
\hspace{.5in}%
\subfigure{\includegraphics[width=1.5in]{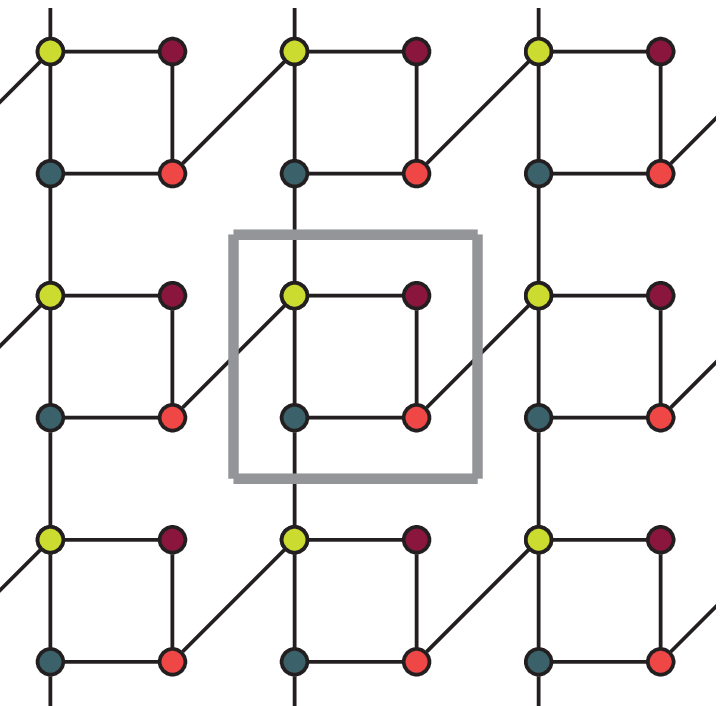}}
\caption{A periodic orbit graph $\pog$, where $m:E \rightarrow \mathbb Z^2$, and the corresponding periodic framework. Any labeled edge in $\pog$ corresponds to an edge in the periodic framework which crosses the boundary of the ``unit cell" (grey box) marked in (b).   \label{fig:gainGraph}}
\end{center}
\end{figure}

A {\it periodic framework} can be described by a locally finite infinite graph $\tilde G$, together with a periodic position of its vertices $\tilde p$ in $\mathbb R^d$ such that the resulting (infinite) framework is invariant under a symmetry group $\Gamma$, which contains as a subgroup the $d$-dimensional integer lattice $\mathbb Z^d$ \cite{BorceaStreinu}. A {\it periodic orbit framework} $\pofw$ consists of a {\it periodic orbit graph} $\pog$ together with a position of its vertices onto the ``flat torus" $\mathcal T^d = \mathbb R^d / \mathbb Z^d$. The periodic orbit graph is a finite graph $G$ which is the quotient of $\tilde G$ under the action of $\Gamma$, together with a labeling of the directed edges of $G$, $m: E(G)^+ \rightarrow \mathbb Z^d$. This periodic orbit framework provides a ``recipe" for the larger periodic framework, but does so with a finite graph $G$, which we can then consider using inductive constructions (Figure \ref{fig:gainGraph}). In addition, it is possible to define a {\it generic} position of the framework vertices on the torus $\mathcal T^d$.

\subsection{Fixed Torus}

\begin{figure}[h!]
\begin{center}
\includegraphics[width=4in]{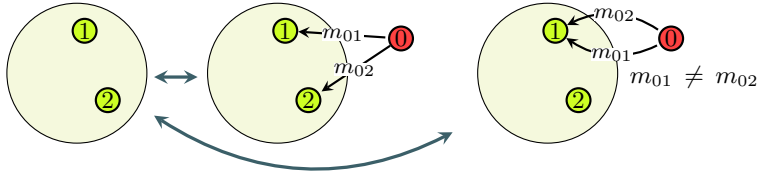}
\caption{Periodic vertex addition. The large circular region represents a generically rigid periodic orbit graph.  \label{fig:vertexAddition}}
\end{center}
\end{figure}

\begin{figure}
\begin{center}
\includegraphics[width=4.5in]{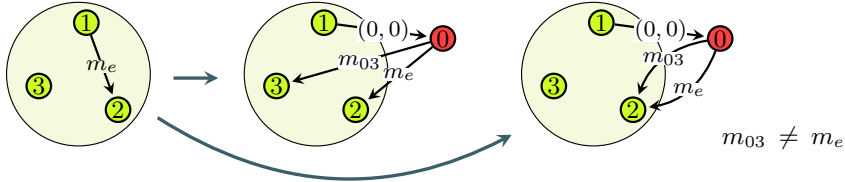}
\caption{Periodic edge split. The net gain on the edge connecting vertices $1$ and $2$ is preserved. \label{fig:edgeSplit}}
\end{center}
\end{figure}

The torus $\T^2$ in $2$ dimensions can be seen as being generated by two lengths and an angle between then. When we do not allow the lengths or angle to change, we call the resulting structure the {\it fixed torus}, and denote it $\Tor^2$. In \cite{RossThesis}, a Laman-type characterisation of graphs which are minimally rigid on the fixed torus is obtained. The proof depended on the development of inductive constructions on periodic orbit graphs $\pog$. These moves require an additional layer of complexity over the usual vertex addition and edge splitting operations. The directed, labeled edges of $\pog$ are recorded by $e = \{v_1, v_2; m_e\}$. We have the following moves: 

\begin{defn}
Let $\pog=(V\pog, E\pog)$ be a periodic orbit graph $\pog$, a {\rm periodic vertex addition} is the addition of a single new vertex $v_0$ to $V\pog$, and the edges $\{v_0, v_{i_1}; \bm_{01}\}$ and $\{v_0, v_{i_2}; \bm_{02}\}$ to $E\pog$, such that $\bm_{01} \neq \bm_{02}$ whenever $v_{i_1} = v_{i_2}$ (see Figure \ref{fig:vertexAddition}). 

Let $e = \{v_{i_1}, v_{i_2}; \bm_e\}$ be an edge of $\pog$. A {\rm periodic edge split} $\pogp$ of $\pog$ results in a graph with vertex set $V \cup \{v_0\}$ and edge set consisting of all of the edges of $E\pog$ except $e$, together with the edges \[ \{v_0, v_{i_1}; (0,0)\}, \{v_0, v_{i_2}; \bm_e\}, \{v_0, v_{i_3}; \bm_{03}\}\]
where $v_{i_1} \neq v_{i_3}$, and $\bm_{03} \neq \bm_{e}$ if $v_{i_2}=v_{i_3}$ (see Figure \ref{fig:edgeSplit}). 
\end{defn}

Together the periodic vertex addition and edge split characterise generic rigidity on the fixed two-dimensional torus $\Tor^2$. Note that the single vertex graph $\pog$ is generically rigid on $\Tor^2$. 

\begin{thm}[Ross \cite{RossThesis}]
A periodic orbit framework $\pofw$ on $\Tor^2$ is generically minimally rigid if and only if it can be constructed from a single vertex on $\Tor^2$ by a sequence of periodic vertex additions and edge splits. 
\label{thm:Henneberg}
\end{thm}

\subsection{Partially Variable Torus}

In \cite{RossThesis}, a characterisation was established of the generic rigidity of periodic frameworks on a partially variable torus (allowing one degree of flexibility). Recently, the authors of the present paper have outlined an inductive proof of this result \cite{NixonRoss}. 

\begin{thm}[Nixon and Ross \cite{NixonRoss}]
A framework $\pofw$ is generically minimally rigid on the partially variable torus (with one degree of freedom) if and only if it can be constructed from a single loop by a sequence of gain-preserving Henneberg operations. 
\label{thm:partFlexHenn}
\end{thm}

The operations referred to in Theorem \ref{thm:partFlexHenn} contain the periodic vertex addition and edge split operations described above. However, we also require one additional move, which is only used in a particular special case. It is an infinite but controllable class of graphs for which vertex addition and edge splitting is insufficient. In addition, while all the generically rigid graphs on the partially variable torus are $(2,1)$-tight, in fact the class of generically rigid graphs is strictly smaller. It is the set of graphs which can be decomposed into an edge-disjoint spanning tree and a {\it connected} spanning map-graph (a connected graph contained exactly one cycle). This hints at the subtlety involved when moving from graphs on the fixed torus to graphs on a partially variable torus, and suggests some challenges which may exist in trying to inductively characterise graphs on the fully flexible torus. 

\subsection{Fully Flexible Torus}

Generic minimal rigidity on the fully variable torus has been completely characterised by Malestein and Theran \cite{MalesteinTheran}. Their proof is non-inductive, however, and there remain significant challenges to providing such a constructive characterisation, since the underlying orbit graph may have minimum degree 4. As we have seen the $X$- and $V$-replacement moves, are known to be problematic in other settings \cite{Whiteleyscene}. 

It may be possible to define somewhat weaker versions of inductive constructions in these settings, by relaxing our focus on ``gain-preservation". That is, we can perform vertex addition and edge splitting on the orbit graph, but allow relabeling of the edges. This is, in some ways, a less satisfying if easier approach, as the moves no longer correspond to the ``classical" inductive moves on the (infinite) periodic framework.

\section{Symmetric Frameworks}
\label{sec:symmetric}

A second class of frameworks which have experienced increased attention over the past decade is symmetric frameworks \cite{RoyalSociety,MalesteinTheranCryst,Schulze,Schulze2}, and there are connections with the study of protein structure. Like periodic frameworks,  symmetric frameworks are frameworks which are invariant under the action of certain symmetry groups, in this case, finite point groups.  

Inductive constructions played a key role in Schulze's work on symmetric frameworks \cite{Schulze}. A {\it symmetric framework} is a finite framework $(G, p)$ which is invariant under some symmetric point group. In $2$-dimensions, this could be for example $\C_2$, half-turn symmetry or $\C_{\S}$, mirror symmetry. Schulze used symmetrized versions of vertex addition and edge splitting to prove Henneberg and Laman-type results for several classes of symmetric frameworks in $\mathbb R^2$, namely $\C_2, \C_3$, and $\C_s$. Furthermore, these results are stronger than the analogous results in the periodic setting, in that they are concerned with frameworks which are either forced to be symmetric, or frameworks which are simply incidentally symmetric. That is, the symmetry-adapted moves preserve the rank of both the (symmetry) orbit matrix, and of the original rigidity matrix of any given symmetric framework.  

\begin{figure}\begin{center}
\includegraphics[width=3.5in]{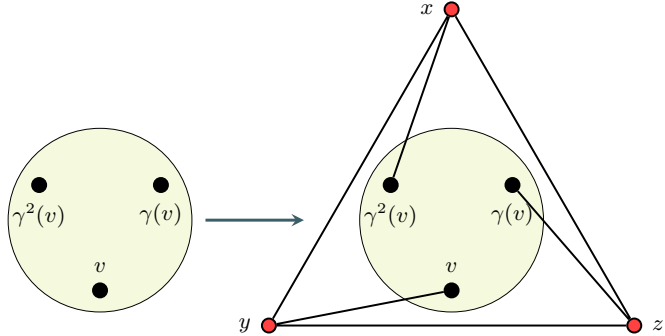}
\caption{One of the three $\C_3$-symmetric edge splitting operations, where $\gamma$ represents rotation through $2\pi/3$. Henneberg proved that the natural generalization of this move preserves rigidity for arbitrary $n$-gons \cite{Henneberg}. It should be noted, however, that Schulze proved the $\C_3$ move for the non-generic ``special geometric" position shown above, where $y = \gamma(x), z = \gamma^2(x)$, and his arguments could easily be extended to cover non-generic $n$-gons (under $\C_n$ symmetry) as well.}
\label{fig:deltaMove}
\end{center}
\end{figure}

As an example, we consider a framework with three-fold rotational symmetry (the group $\C_3$). 
\begin{thm}[Schulze \cite{Schulze}]
A $C_3$-symmetric framework $(G,p)$ is generically (symmetric)-isostatic if and only if it can be generated through three inductive moves, a three-fold vertex addition (one vertex is added symmetrically to each of the three orbits), a three-fold edge split (one edge is ``split" symmetrically in each of the three orbits) and the $\Delta$-move pictured in Figure \ref{fig:deltaMove}. 
\end{thm}
Schulze proves analogous results for $\C_2$ and $\C_s$ \cite{Schulze2}. In the case of $\C_s$ (mirror symmetry), $X$-replacement is also required to handle certain special cases. Schulze also proves tree-covering results for these groups. 

We remark that it would be possible to rework these results of Schulze using the language of {\it gain graphs} (graphs whose edges are labeled by group elements), as for periodic frameworks. In that scenario, we would capture the symmetric graph using an orbit graph whose edges were labeled with elements of the symmetry group (e.g. $\C_3$ etc.). The symmetric inductive moves could then be defined on this symmetric orbit graph. This is exactly the approach taken in very recent work of Jord{\'a}n, Kaszanitsky and Tanigawa \cite{JordanKaszanitskyTanigawa}, for the groups $\C_s$ (the reflection group), and the dihedral groups $D_h$, where $h$ is odd. We mention here their results for $D_h$. 

The authors define a $D_h$ sparsity type of the gain graphs $(G, \phi)$, where $\phi$ is a labeling of the edges by elements of the group $D_h$. They then prove that all $D_h$-tight graphs can be constructed from the disjoint union of a few `basic' graphs by a sequence of Henneberg-type moves on the underlying gain graph. In particular, they use vertex addition, edge splitting and $X$-replacements; including {\it loop vertex addition} (adding a `lolipop'), and {\it edge splitting plus adding a loop on the new vertex}. 
This leads to the following combinatorial characterisation of rigid frameworks with $D_h$ symmetry (A similar result is established for $\C_s$):

\begin{thm}[Jord\'{a}n, Kaszanitsky and Tanigawa \cite{JordanKaszanitskyTanigawa}]
$(G, \phi)$, $\phi: E(G) \rightarrow D_h$, where $h$ is odd, is the gain graph of a rigid framework with $D_h$ symmetry if and only if $(G, \phi)$ has a $D_h$-tight subgraph. 
\end{thm}
 
Note that the work of Schulze provides combinatorial characterisations for $\C_2$, $\C_3$, and $\C_s$ only, but his results are for both incidental and forced symmetry. On the other hand, Jord\'{a}n, Kaszanitsky and Tanigawa's results are for forced symmetry only but cover any cyclic group and odd order dihedral groups. Thus, there are a number of outstanding questions about symmetric frameworks, including the characterisation of the rigidity of frameworks with forced dihedral $D_h$ ($h$ even) symmetry, and the characterisation of incidental rigidity for the dihedral groups. 

We remark that for the cyclic groups $C_n$ and $C_s$ (rotations and reflections) the characterisations of forced rigidity can also be obtained using direction networks and linear representability \cite{MalTheran}.

\section{Frameworks on Surfaces}
\label{sec:surfaces}

Inductive constructions have also played a big role in recent work on frameworks supported on surfaces. Here characterisations of minimal rigidity require us to remain within the class of simple graphs at each step of the induction. Hence the results of Fekete and Szeg\H{o}, allowing multiple edges and loops, are not sufficient. For example a $(2,2)$-tight graph may contain an arbitrarily large number of copies of $K_4$ and there is no inverse edge splitting operation on a degree $3$ vertex in a copy of $K_4$ that preserves simplicity.

\begin{center}
\begin{figure}[ht]
\centering
\includegraphics[width=4.2cm]{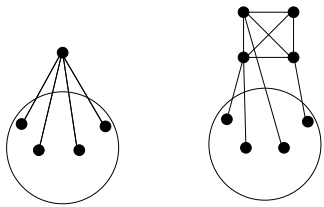}
\caption{An example of the vertex-to-$K_4$ operation.}
\label{fig:vtok4}
\end{figure}
\end{center}

This motivates the vertex-to-$K_4$ move, in which we remove a vertex $v$ (of any degree) and all incident edges $vx_1,\dots, vx_n$ and insert a copy of $K_4$ along with edges $x_1y_1,\dots,x_ny_n$ where each $y_i \in V(K_4)$, see Figure \ref{fig:vtok4}.

\begin{thm}[Nixon and Owen \cite{NixonOwen}]
A simple graph $G$ is $(2,2)$-tight if and only if $G$ can be generated from $K_4$ by vertex addition, edge splitting, vertex-to-$K_4$ and ($2$-dimensional) vertex splitting operations.
\end{thm}

Similarly when dealing with $(2,1)$-tight graphs, all low degree vertices may be contained in copies of $K_5-e$ (the graph formed from $K_5$ by deleting any single edge). For these graphs vertex-to-$K_4$ and vertex splitting moves are not sufficient so we introduce the edge joining move. This is the joining of two $(2,1)$-tight graphs by a single edge.
In the following theorem $K_4 \sqcup K_4$ is the unique graph formed from two copies of $K_4$ intersecting in a single edge.

\vspace{-1.2cm}

\begin{thm}[Nixon and Owen \cite{NixonOwen}]
A simple graph $G$ is $(2,1)$-tight if and only if $G$ can be generated from $K_5-e$ or $K_4 \sqcup K_4$ by vertex addition, edge splitting, vertex-to-$K_4$, vertex splitting and edge joining operations.
\end{thm}

The first of these results has led to an analogue of Laman's theorem for an infinite circular cylinder and the second to analogues for surfaces admitting a single rotational isometry (such as the cone and torus).

\begin{thm}[Nixon, Owen and Power \cite{NixonOwenPower}]
Let $G=(V,E)$ with $|V|\geq 4$. Then the framework $(G,p)$ is generically minimally rigid on a cylinder if and only if $G$ is simple and $(2,2)$-tight.
\end{thm}

\begin{thm}[Nixon, Owen and Power \cite{NixonOwenPower2}]
Let $G=(V,E)$ with $|V|\geq 5$. Then the framework $(G,p)$ is generically minimally rigid on a surface of revolution if and only if $G$ is simple and $(2,1)$-tight.
\end{thm}

The next extension would be to frameworks on a surface admitting no ambient isometries (such as an ellipsoid). This would require an inductive construction of $(2,0)$-tight graphs; these have average degree $4$, in fact they may be $4$-regular so $X$- and $V$- replacement and their associated complication may be required.

The insistence on simplicity also makes the characterisation of $(k,l)$-circuits more challenging. Similarly to Berg and Jord{\'a}n's theorem the following inductive result is a step towards characterising global rigidity on the cylinder. The $1$-, $2$- and $3$-join operations, defined in \cite{Nixon}, are similar in spirit to the 2-sum operation used in Theorem \ref{thm:BergJordancircuit}. 

\begin{thm}[Nixon \cite{Nixon}]
A simple graph $G$ is a $(2,2)$-circuit if and only if $G$ can be generated from copies of $K_5-e$ and $K_4 \sqcup K_4$ by applying edge splitting within connected components and taking $1$-, $2$- and $3$-joins of connected components.
\end{thm}

It is an open problem to extend this characterisation to give an inductive construction for generically globally rigid frameworks on a cylinder.

\section{Body-bar Frameworks}
\label{sec:bodyBar}

Body-bar frameworks are a special class of frameworks where there is a more complete understanding in arbitrary dimension. Roughly speaking, a body-bar framework is a set of bodies (each spanning an affine space of dimension at least $d-1$), which are linked together by stiff bars. 
\begin{thm}[Tay \cite{TayOriginal}]
Let $G$ be a graph. Then $(G,p)$ is generically minimally rigid as a body-bar framework in $\mathbb{R}^d$ if and only if $G$ is $(D,D)$-tight, where $D = \binom{d+1}{2}$ is the dimension of the Euclidean group. 
\end{thm}

Tay subsequently proved an inductive characterisation of the body-bar frameworks. 

\begin{thm}[Tay \cite{Tayhenneberg}]
A graph $G$ is $(D,D)$-tight if and only if $G$ can be formed from $K_1$ by Henneberg operations.
\label{thm:tayHenn}
\end{thm}
The Henneberg operations referred to in Theorem \ref{thm:tayHenn} are essentially the loopless versions of the edge-pinches of Fekete and Szeg\H{o} (see Section \ref{sec:klTight}).  

Recently, Katoh and Tanigawa proved the {\it Molecular Conjecture}, a longstanding open question due to Tay and Whiteley, which is concerned with body-bar frameworks which are geometrically special:
\begin{thm}[Katoh and Tanigawa \cite{KatohTanigawa}] Let $G = (V, E)$ be a graph. Then, $G$ can be realised as an infinitesimally rigid body-and-hinge framework in $\mathbb{R}^d$ if and only if $G$ can be realised as an infinitesimally rigid panel-and-hinge framework in $\mathbb{R}^d$.
\end{thm}
The settling of this conjecture is of particular significance to the materials science community, who use rigidity analysis for the modeling of molecular compounds. The proof of this result is quite involved, so we will not include many details here. However, one of the ingredients in the proof is inductive constructions. In particular, the authors use a type of {\it splitting off} operation, which removes a two-valent vertex $v$, and then inserts a new edge between the pair of vertices formerly adjacent to $v$. A second type of induction used is a {\it contraction} operation, which contracts a proper rigid subgraph to a vertex. 

Along the way Katoh and Tanigawa also obtain a Henneberg-type characterisation of minimally rigid body-and-hinge graphs. In particular, they show that for any minimally rigid body-and-hinge framework, there is a sequence of graphs ending with the two vertex, two edge graph, where each graph in the sequence is obtained from the previous graph by a splitting off operation or a contraction operation (see Theorem 5.9, \cite{KatohTanigawa}). 

\subsection{Global Rigidity}

Inductive constructions have also played a role in the proof of the following result concerning generic global rigidity of body-bar frameworks:
\begin{thm}[Connelly, Jord{\'a}n and Whiteley \cite{genericGlobal}] 
A body-bar framework is generically globally rigid in $\bR^d$ if and only if it is generically redundantly rigid in $\bR^d$.
\label{thm:genericGlobal}
\end{thm}
In particular, the authors' proof used Theorem \ref{thm:frankszego} to produce an inductive construction of redundantly rigid body-bar graphs. One of the interesting elements of this proof is that the construction sequence specified by Theorem \ref{thm:frankszego} may involve loops. However, no (globally) rigid finite framework will involve loops. The proof of Theorem 
\ref{thm:genericGlobal}  involved allowing for the possibility of loops, which would later be eliminated. In this way, the induction used here stepped outside of the class of frameworks under study, but eventually achieved the desired result. 

\subsection{Periodic Body-Bar Frameworks}
\label{sec:rossBodyBar}

It is possible to define periodic body-bar frameworks in much the same way as periodic bar-joint frameworks (see Section \ref{sec:periodic}). A recent result of Ross characterises the generic rigidity of periodic body-bar frameworks on a three dimensional fixed torus \cite{RossBodyBar}. It is based on the following sparsity condition which depends on the dimension of the {\it gain space} $\gs$: the vector space generated by the net gains on all of the cycles of a particular edge set $Y$. 
\begin{thm}[Ross \cite{RossBodyBar}]
$\bbog$ is a periodic orbit graph corresponding to a generically minimally rigid body-bar periodic framework in $\mathbb R^3$ if and only if  $|E(H)| = 6|V(H)| - 3$ and for all non-empty subsets $Y \subset E(H)$ of edges
		\begin{equation*}
			|Y| \leq 6|V(Y)| - 6 + \sum_{i=1}^{|\gs(Y)|}(3-i).
		\end{equation*}	
\end{thm}
The proof relies on a careful modification of the edge-pinching results of Fekete and Szeg\H{o} \cite{FeketeSzego} to include labels on the edges of the graphs. It is interesting to note that the results of Fekete and Szeg\H{o} cover the class of minimally rigid frameworks on the fixed torus, but will not assist us with the flexible torus. That is, for minimal rigidity on the fixed torus, we are considering $(\binom{d+1}{2}, d)$-tight graphs, whereas for minimal rigidity on the flexible torus, we are considering $(\binom{d+1}{2}, -\binom{d}{2})$-tight graphs, which are not in the range covered by existing inductive results. Periodic body-bar frameworks with a flexible lattice have recently been considered in \cite{PeriodicBarBody} using non-inductive methods. 

\section{Further Inductive Problems}
\label{sec:further}

Aside from the conjectures already discussed, a number of other problems, especially in $3$-dimensions, remain open, see \cite{Taygeneric}, \cite{GraverServatius2}, \cite{Whiteleymatroid}. 

There are a number of connections between two-dimensional minimally rigid frameworks and the topic of pseudo-triangulations. A {\it pseudo-triangulation} is a tiling of a planar region into {\it psuedo-triangles}: simple polygons in the plane with exactly three convex vertices \cite{RoteSantosStreinu}. It is called a {\it pointed} pseudo-triangulation if every vertex is incident to an angle larger than $\pi$. Streinu proved that the underlying graph of a pointed pseudo-triangulation of a point set is minimally rigid \cite{ParallelRedrawing}. As a converse, there is the following result: 
\begin{thm}[Haas et. al. \cite{PseudoTriangles}]
Every planar infinitesimally rigid graph can be embedded as a pseudo-triangulation. 
\end{thm}
The proof uses vertex addition and edge splitting. Pseudo-triangulations are the topic of an extensive survey article \cite{RoteSantosStreinu}, and further details on the inductive elements of the proof can be found there. 

In \cite{PilaudSantos} Pilaud and Santos consider an interesting application of rigidity in even dimensions to multitriangulations. In particular they use Theorem \ref{thm:vertexsplittingd} to show that every $2$-triangulation is generically minimally rigid in $4$-dimensions and conjecture the analogue for $k$-triangulations in $2k$-dimensions.

If a framework is not globally rigid then the number of equivalent realisations of the graph is not unique. For $d\geq 2$ this is not a generic property, nevertheless bounds on the number of realisations were established by Borcea and Streinu \cite{BorceaStreinurealisations} and recent work of Jackson and Owen \cite{JacksonOwen}, motivated by applications to Computer Aided Design (CAD), considering the number of complex realisations made use of vertex addition and edge splitting.

Servatius and Whiteley \cite{ServatiusWhiteley}, again motivated by CAD, used the Henneberg operations to understand the rigidity of direction-length frameworks. Jackson and Jord{\'a}n \cite{JacksonJordandirection} established the analogue of Theorem \ref{thm:BergJordancircuit} for direction-length frameworks however a characterisation of globally rigid direction-length frameworks remains an open problem.

\bibliographystyle{abbrv} 
\bibliography{inductivesurveyRevisedVersion}
\end{document}